\documentclass[12pt]{article}
\usepackage{ amsfonts,
amsmath, amssymb,amsgen, amsthm, amscd}

\def\GL{\mathop{\rm GL}\nolimits}

\def\Hopf{\mathop{\rm Hopf}\nolimits}
\def\Id{\mathop{\rm Id}\nolimits}

\def\Ker{\mathop{\rm Ker}\nolimits}

\def\log{\mathop{\rm log}\nolimits}

\def\SL{\mathop{\rm SL}\nolimits}
\def\PSL{\mathop{\rm PSL}\nolimits}

\def\Ac{{\cal A}}

\def\Fc{{\cal F}}

\def\Hc{{\cal H}}

\def\Ac{{\cal A}}

\def\Fc{{\cal F}}

\def\Hc{{\cal H}}

\def\Lc{{\cal L}}
\def\Mc{{\cal M}}
\def\Pc{{\cal P}}

\def\a{\alpha}
\def\b{\beta}
\def\d{\delta}

\def\D{\Delta}
\def\g{\gamma}

\def\om{\omega}
\def\Om{\Omega}
\def\s{\sigma}

\def\ve{\varepsilon}

\def\tb{\bf t}

\def\0b{\bf 0}

\def\fl{\forall}

\def\ot{\otimes}

\def\ra{\rightarrow}

\def\wt{\widetilde}

\newtheorem{theorem}{Theorem}

\newtheorem{proposition}[theorem]{Proposition}
\newtheorem{lemma}[theorem]{Lemma}
\newtheorem{corollary}[theorem]{Corollary}

\def\build#1_#2^#3{\mathrel{ 
\mathop{\kern 0pt#1}\limits_{#2}^{#3}}}

\numberwithin{equation}{section}
\begin{document}
\title{{\bf  Rankin-Cohen Brackets and  the Hopf
Algebra of Transverse Geometry}}
\author{Alain Connes \\ 
        Coll\`ege de France \\
        3 rue d'Ulm \\ 
        75005 Paris, France \\
\and 
        Henri Moscovici\thanks{Research 
    supported by the National Science Foundation
    award no. DMS-0245481.} \\
    Department of Mathematics \\
    The Ohio State University \\
    Columbus, OH 43210, USA
}
 
\date{ \ }

\maketitle

\centerline{\it Dedicated to Pierre Cartier}

\begin{abstract}
\noindent We settle in this paper a question
left open in  our paper ``Modular Hecke algebras
and their Hopf symmetry'', by showing how to extend the 
Rankin-Cohen brackets from modular forms to modular Hecke algebras.
More generally, our procedure yields such brackets on any
associative algebra 
endowed with an action of the Hopf algebra of transverse geometry
in codimension one, 
such that the derivation corresponding to the Schwarzian derivative is inner.
Moreover, we show in full generality that these Rankin-Cohen brackets give rise to
associative deformations.
\end{abstract}

\section*{Introduction}

\noindent In~\cite{HC} H. Cohen constructed a collection of
bilinear operations on functions on the complex plane,
which are covariant with respect to the `slash' action of
$\PSL(2, \mathbb{R})$ and therefore give rise to 
bi-differential operators on the graded algebra of modular forms.
If $f$ and $g$ denote two modular forms of weight $k$ and $\ell$,
the $n$th \textit{Rankin-Cohen bracket} of $f$ and $g$ is 
given (using the normalization in~\cite{Za}) by the formula
\begin{equation}
[f, g]_n \, := \, \sum_{r+s=n} (-1)^r \binom{n+k-1}{s}
\binom{n+\ell-1}{r} f^{(r)} g^{(s)} \, .
\end{equation}
Don Zagier later investigated
the abstract algebraic structure defined by such
bilinear operations~\cite{Za}.

\medskip

\noindent  We noticed in~\cite{CM10} that
the formulas for the perturbations by inner of the
action of the Hopf algebra $\Hc_1$ of \cite{CM1}
on the modular Hecke algebras were similar to the 
formulas occuring in Zagier's definition~\cite{Za} of 
`canonical' Rankin-Cohen algebras.
This suggested that there should be 
a close relationship between Rankin-Cohen
brackets and actions of the Hopf algebra $\Hc_1$.
We shall show here that this is indeed the case,
by exhibiting canonical bilinear operators that
generalize the Rankin-Cohen brackets,
for any associative algebra $\Ac$ endowed with 
an action of the Hopf algebra $\Hc_1$ of \cite{CM1}
for which the `Schwarzian' derivation $\d'_2$ 
of \cite{CM10} is inner,
\begin{equation} \label{inner}
\d'_2(a)= \Omega \, a - \, a\, \Omega\, ,\qquad \forall
a \in \Ac \, .
\end{equation}
Such an action of the Hopf algebra $\Hc_1$
defines the noncommutative analogue of a one-dimensional
``projective structure", while $\Omega$ plays the 
role of the quadratic differential.

\noindent We had given in  \cite{CM10} the first term in the 
deformation as (twice) the fundamental class $[F]$ in the cyclic 
cohomology $PHC_{\Hopf}^{\rm ev}(\Hc_1)$,
i.e. the class of the cyclic $2$-cocycle
\begin{equation} \label{fund}
F := X \ot Y - Y \ot X - \d_1 \,  Y \ot Y \, ,
\end{equation}
which in the foliation context
represents the transverse fundamental class.
Since the antipode $S(X)$ is given by 
\begin{equation} \label{ant}
S(X)= -X + \, \d_1 \: Y  \, ,
\end{equation}
one has 
$$ -F \,=\,  S(X) \ot Y + Y \ot X \, ,
$$ 
and one could reasonably expect the higher brackets to 
be increasingly more complicated expressions
involving $\displaystyle S(X), \, X, \, Y \,$ and $\Omega$,
beginning with the first bracket
\begin{equation} \label{rc}
RC_1(a,b):= S(X)(a) \, 2Y(b)\,+ 2Y(a) \, X(b) \, ,
\qquad a , b \in \Ac \, .
\end{equation}

\noindent We shall obtain in this paper the  
canonical formulas for the higher brackets. 
 As an illustration, the 
second bracket $RC_2$ is given by
\begin{eqnarray}
 RC_2(a,b):=&&S(X)^2(a) \, Y ( 2\, Y+1 )(b)\,+
\,S(X)\,( 2\, Y+1 )(a)\, X ( 2\, Y+1 )(b) \nonumber \\
&& + \,
 Y ( 2\, Y+1 )(a)  \, X^2(b) \, - \,
Y(a)\, \Omega \, Y(2Y+1)(b) \nonumber \\
&& - \, Y(2Y+1)(a) \, \Omega \, Y(b)  
\end{eqnarray}
The ``quadratic differential'' $\Omega$ and 
its higher derivatives $X^j(\Omega)$ always intervene
between polynomial expressions $P(S(X), Y)(a)$ and $Q(X,Y)(b)$. 
We shall also prove (Theorem \ref{cross})
that when applied to the modular Hecke algebras
they yield a family of associative formal deformations, which 
in particular incorporate the 
`tangent groupoid' deformation.
This will be achieved
by showing that the canonical formula commutes with
the crossed product construction under conditions
which are realized for the action of Hecke operators on
modular forms, and then relying on the results of \cite{CMZ}. 
The commutation with crossed product will in fact
uniquely dictate the general formula.
We shall show that the general formula is invariant
under inner perturbations of the action of the 
Hopf algebra.
In order to obtain the whole sequence of  higher components
we shall closely follow the method of \cite{{Za}} and \cite{{CMZ}}. 
\medskip

\noindent Passing to the general algebraic context,
in the last section we establish (Theorem \ref{assoc1})
the associativity of the formal
deformations corresponding to the Rankin-Cohen brackets
for an arbitrary associative algebra 
endowed with an action of $\Hc_1$
such that the derivation corresponding to the Schwarzian derivative is inner.
\medskip

\tableofcontents

\section{Zero Quadratic Differential} 

For the clarity of the exposition, we shall
 first develop the higher Rankin-Cohen brackets 
for actions of $\Hc_1$
in the simplified case when the quadratic differential 
$\Omega$ is zero, i.e.  when $\d'_2=0$.
\smallskip

\noindent We recall that
as an algebra $\Hc_1$ coincides with the universal 
enveloping algebra of the Lie algebra with basis $\{ X,Y,\d_n \, ; n \geq 1 
\}$ and brackets
\begin{equation*}
[Y,X] = X \, , \, [Y , \d_n ] = n \, \d_n \, , \, [X,\d_n] = \d_{n+1} 
\, , \,
[\d_k , \d_{\ell}] = 0 \, , \quad n , k , \ell \geq 1 \, ,
\end{equation*}
while the coproduct which confers it the Hopf algebra structure 
is determined by the identities
\begin{eqnarray}
\D \,  Y &=& Y \ot 1 + 1 \ot Y \, , \quad
\D \,  \d_1 = \d_1 \ot 1 + 1 \ot \d_1 \, , \nonumber
 \\
\D \,  X &=& X \ot 1 + 1 \ot X + \d_1 \ot Y \, , \nonumber
\end{eqnarray}
together with the property that 
$\, \D : \Hc_1 \ra \Hc_1 \ot \Hc_1 \, $
is an algebra homomorphism.

\noindent We let 
 $\Hc_{s}$ be the  quotient of $\Hc_1$ by 
the ideal generated by $\delta'_2$.
The basic lemma  is to put the antipode
$S(X^n)$ in normal order, i.e. with the 
$\delta's$ first then the $X$ then the $Y's$.

\begin{lemma} \label{start} One has the following identity
in $\Hc_{\rm s}$,
$$
S(X^n) = \sum (-1)^{n-k}\,\binom{n}{k} \, \,\frac{\delta_1^k}{2^k} \,
X^{n-k} (2Y+n-k)_k
$$
where $\displaystyle (\alpha)_k := \alpha \,(\alpha+1)....(\alpha+k-1)$.
\end{lemma}
\begin{proof} For $n=1$ one has $\displaystyle S(X)= -X + \d_1 \, Y$ 
and the r.h.s. of the formula gives 
$\displaystyle -X + \,\frac{\delta_1}{2} \, (2Y)$.
Let us assume that it holds for $n$ and check it for $n+1$
by multiplication on the left by $\displaystyle S(X)= -X + \d_1 \, Y$. 
One has $\displaystyle Y \, \frac{\delta_1^k}{2^k} \, X^{n-k}=\, 
\frac{\delta_1^k}{2^k} \, X^{n-k} (Y+n)$
so that the new terms corresponding to $\displaystyle 
\d_1 \, Y \, S(X^n) $ give
$$
\sum (-1)^{n-k}\,\binom{n}{k} \, \,\frac{\delta_1^{k+1}}{2^{k+1}} \, 
X^{n-k} (2Y+2n)(2Y+n-k)_k
$$
One has $\displaystyle [X, \, \delta_1^k]=\, k \,\frac{\delta_1^{k+1}}{2}$ 
thus the new terms corresponding to $\displaystyle 
-X \, S(X^n) $ give
$$
\sum (-1)^{n+1-k}\,\binom{n}{k} \, \,\frac{\delta_1^k}{2^k} \, 
X^{n+1-k} (2Y+n-k)_k
$$
and
$$
\sum (-1)^{n+1-k}\,\binom{n}{k} \, \,\frac{\delta_1^{k+1}}{2^{k+1}} \, 
X^{n-k} \, k\,(2Y+n-k)_k
$$
The first of these two expressions can be written as
$$
\sum (-1)^{n-k}\,\binom{n}{k+1}  \, \,\frac{\delta_1^{k+1}}{2^{k+1}} \, 
X^{n-k} (2Y+n-k-1)_{k+1}
$$
thus all the term are right multiples of
$$
(-1)^{n-k}\,\frac{\delta_1^{k+1}}{2^{k+1}} \, X^{n-k} \,
(2Y+n-k)_k
$$
with coefficients
$$
\,\binom{n}{k}\, (2Y+2n)\, -\,\binom{n}{k}\,k +\,\binom{n}{k+1}\,(2Y+n-k-1)
$$
but this is the same as 
\begin{eqnarray}
&& \left(\,\binom{n}{k}\,+\,\binom{n}{k+1}\,\right)\,(2Y+n)\,+ 
\, (n-k)\,\binom{n}{k}\,
-\,(k+1)\, \binom{n}{k+1}  \nonumber \\
&& =\,\binom{n+1}{k+1}\, \,(2Y+n)\, \nonumber
\end{eqnarray}

\end{proof}
\noindent The content of Lemma \ref{start} can be written using the 
generating function 
\begin{equation}
\Phi(X)(s):=  \sum \frac{s^n \, X^n}{n!}\, \Gamma(2Y+n)^{-1}
\end{equation}
as the equality,
\begin{equation} \label{pert1}
\Phi(X -Z \, Y)(s)=  e^{-\frac{s\, Z}{2}}\,\Phi(X)(s)
\end{equation}
under the only assumption that the operators $X,Y,Z$
fulfill,
\begin{equation}\label{pert0}
[Y,X]=X \,, \quad [Y,Z]=Z \,, \quad [X,Z]=\frac{1}{2}Z^2 \,, \end{equation}

\noindent Expanding the tensor product
\begin{equation}\label{rc1}
\Phi(S(X))(s)\ot  \Phi(X)(s)=\sum s^n \,
RC^{\Hc_{\rm s}}_n \:\;(\Gamma(2Y+n)\ot \Gamma(2Y+n))^{-1}
\end{equation}
 dictates the following formula for
the higher Rankin-Cohen brackets for an action of $\Hc_{s}$ on an algebra $\Ac$,
\begin{equation} \label{RC0}
RC_n(x,\,y):=\, \sum_{k=0}^n \, (\frac{S(X)^k}{k!} \,
(2Y+k)_{n-k})(x)(\frac{X^{n-k}}{(n-k)!} \,(2Y+n-k)_{k})(y)   
\end{equation}

\noindent The main preliminary result then is the following.

\begin{lemma} \label{11} Let the Hopf algebra $\Hc_1$ act on an 
algebra $\Ac$ and $u \in \Ac$ be invertible and such that,
$$
X(u)=0, \, Y(u)=0, \, \d_1(u^{-1}\d_1(u))=0, \, \d'_2(u)=0\,.
$$
\begin{itemize}
\item[$1^0$.] For all values of $n$,
$$
RC_n(x\,u,\,y)=\, RC_n(x,\,u\,y) \quad\forall x,\, y \in \Ac
$$

\item[$2^0$.] For all $x,\, y \in \Ac
$,
$$
RC_n(u\,x,\,y)=\,u\, RC_n(x,\,\,y)\, ,
\quad RC_n(x,\,y \,u)=\, RC_n(x,\,\,y)\, u  
$$

\item[$3^0$.] Let $\a$ be the inner automorphism implemented by $u$ one has,
$$
RC_n(\a(x),\,\a(y))=\, \a(RC_n(x,\,y)) \quad\forall x,\, y \in \Ac
$$
\end{itemize}
\end{lemma}

\begin{proof}
For $x \in \Ac $ let $L_x$ be the operator of left multiplication by
$x$ in $\Ac $. One has
\begin{equation} \label{left1}
X\,L_u = L_u (\, X - L_{\nu}\, Y)
\end{equation}
where 
 $\nu :=-u^{-1}\d_1(u)= \d_1(u^{-1}) u$. Moreover
$$
Y( \nu )=\nu \, , \quad X(\nu )= \frac{1}{2} \nu^2\,.
$$
so that the operators $X,Y,Z:=L_{\nu}$ fulfill (\ref{pert0})
and by (\ref{pert1})
\begin{equation} \label{pert2}
\Phi(X -L_{\nu} \, Y)(s)= L_{ e^{-\frac{s\, \nu}{2}}}\,\Phi(X)(s)
\end{equation}
One has
$$
\Delta(S(X))= S(X)\ot 1 + 1\ot S(X) + Y \ot \d_1
$$
and using $R_x$ for right multiplication by $x \in \Ac$,
\begin{equation} \label{right1}
S(X)\,R_u = R_u (\, S(X) - R_{\nu'}\, Y)
\end{equation}
where 
 $\nu' :=u\d_1(u^{-1})$. Moreover $\displaystyle \d_1(\nu')=0$
and
$$
Y( \nu' )=\nu' \, , \quad S(X)(\nu' )= \frac{1}{2} (\nu')^2\,,
$$
so that $\, S(X), \, Y, \, R_{\nu'}\,$ fulfill
(\ref{pert0}) and by (\ref{pert1})
\begin{equation} \label{pert4}
\Phi(S(X) -R_{\nu'} \, Y)(s)= R_{ e^{-\frac{s\, \nu'}{2}}}\,\Phi(S(X))(s)
\end{equation}
Since $\displaystyle u\, \nu =\nu' \, u$ we get $1^0$
using (\ref{rc1}) and the commutation of $Y$ with $L_u$ and
$R_u$.
Statement $2^0$ follows from the commutation of 
$R_u$ 
with $X$ and $Y$, while $3^0$ follows from $1^0$ and $2^0$.
\end{proof}

\noindent Let us compute the above brackets for small values of $n$. 
We shall express them as elements 
$$RC^{\Hc_{\rm s}}_n \in \Hc_{\rm s} \ot \Hc_{\rm s} \, .
$$
For $n=1$ we get 
$$
RC^{\Hc_{\rm s}}_1 \,=\,  S(X) \ot 2Y + 2Y \ot X
$$
which is $-2$ times the transverse fundamental
class $F$ (defined in (\ref{fund})).

\noindent
For $n=2$ we get 
\begin{eqnarray}
&& RC^{\Hc_{\rm s}}_2=\, \frac{1}{2} S(X)^2 \ot( 2\, Y )( 2\, Y+1 )\,+
\,S(X)\,( 2\, Y+1 )\ot X \,( 2\, Y+1 )\nonumber \\
&& +
\frac{1}{2}( 2\, Y )( 2\, Y+1 )\, \ot \, X^2 \, . \nonumber
\end{eqnarray}
For $n=3$ we get 
\begin{eqnarray}
&& RC^{\Hc_{\rm s}}_3=\, \frac{1}{6} S(X)^3 
\ot( 2\, Y )( 2\, Y+1 )( 2\, Y+2 )\, \nonumber \\
&& +  \, \frac{1}{2} S(X)^2 \,( 2\, Y+2 )\,  \ot \,
 X\,( 2\, Y+1 )( 2\, Y+2 )\,\nonumber \\
&& +  \, \frac{1}{2} S(X) \,( 2\, Y+1 )\,( 2\, Y+2 )\,  
\ot \, X^2\,( 2\, Y+2 )\, \nonumber \\
&& + 
\,
\frac{1}{6}( 2\, Y )( 2\, Y+1 )\,( 2\, Y+2 )\, \ot \, X^3 \, .
\nonumber
\end{eqnarray}

\section{Arbitrary Quadratic Differential}
\noindent Let us now pass to the general case where
we assume that $\Hc_1$ acts on an algebra $\Ac$ and that
the derivation $\d'_2$ is inner implemented by 
an element $\Omega \in \Ac$ so that,
\begin{equation} \label{om}
\d'_2(a)= \Omega \, a - a \, \Omega\,,\:\forall a \in \Ac
\end{equation}
where we assume, owing to the commutativity of the $\d_k$,
\begin{equation} \label{om1}
   \d_k(\Omega )=\,0\,,\: \forall k \in \mathbb{N}
\end{equation}
It follows then that
\begin{equation} \label{om0}
   \d_k(X^j(\Omega) )=\,0\,,\: \forall k\,, \;j \in \mathbb{N}
\end{equation}
so that by (\ref{om}), $\Omega$ commutes with all $X^j(\Omega)$
and the algebra $\Pc \subset \Ac$
generated by the $X^j(\Omega)$ is commutative,
while both $Y$ and $X$ act as derivations on $\Pc$.

\noindent To understand how to obtain the general formulas
we begin by computing in the above case ($\Omega=0$) how the 
formulas for $RC_n$ get modified by a perturbation 
of the action of the form
\begin{equation}
Y \rightarrow Y \,, \quad X \rightarrow X + \mu \, Y \,,\quad
\d_1 \rightarrow \d_1 + ad(\mu) \,
\end{equation}
where $Y(\mu)=\mu$ and $\d_n(\mu)=0$ for all $n$.
The computation shows that $RC_1$ is unchanged, while $RC_2(a,b)$ 
gets modified
by  the following term,
\begin{equation}
\delta RC_2(a,b)= Y(a)\, \Omega \, Y(2Y+1)(b) + Y(2Y+1)(a)
\, \Omega \, Y(b)
\end{equation}
where $\displaystyle \Omega := X(\mu) + \frac{1}{2} \mu^2$.
Note that the perturbed action fulfills (\ref{om}) 
for that value of $\Omega$.

\noindent This already indicates that in the general case 
the full formula for $RC_2(a,b)$ should be 
\begin{eqnarray}
RC_2(a,b):= &&S(X)^2(a) \, Y ( 2\, Y+1 )(b)\,+
\,S(X)\,( 2\, Y+1 )(a)\, X ( 2\, Y+1 )(b)\nonumber \\
&& + \,
 Y ( 2\, Y+1 )(a)  \, X^2(b) \, - \,
Y(a)\, \Omega \, Y(2Y+1)(b) \nonumber \\
&& - \, Y(2Y+1)(a) \, \Omega \, Y(b) \nonumber 
\end{eqnarray}
so that the above perturbation  then leaves $RC_2$
unaffected.

\noindent In order to obtain the general formulas,
we consider the  algebra $\Lc(\Ac )$
of linear operators in $\Ac$. For $a \in \Ac$ we use the short hand notation 
$$
a :=L_a\,, \quad a^{\,o} :=R_{a}
$$
for the operators of left and right multiplication by 
$a$ when no confusion can arise.

\noindent We define by induction elements  $B_n \in \Lc(\Ac )$ by the equation  
\begin{equation} \label{indb1}
B_{n+1}:= X\, B_n
- n\,\Omega \,(Y- \frac{n-1}{2}) B_{n-1} 
\end{equation}
while $B_0:=1$ and $B_1:=  X$.
The first values for the 
$B_n$ are the following,
$$
B_2=  X^2 - \,\Omega \, Y 
$$
$$
B_3=X^3 - \,\Omega \,X \,(3 Y +1)- \,X(\Omega) \, Y 
$$
and 
$$
B_4=X^4 
-\,\Omega \,X^2 \, (6Y+4) -\,X(\Omega) \,X \, (4Y+1)
- \,X^2(\Omega) \, Y + 3\,\Omega^2 \,Y \,( Y+1)
$$
Let us define more generally for any two operators $Z$ and $\Theta$ acting
linearly in $\Ac$ and fulfilling
$$
[Y,\, Z]=Z\, ,\quad [Y,\, \Theta]=2\Theta
$$
the sequence of operators, $C_0:=1$, $C_1:=Z$,
\begin{equation} \label{indb11}
C_{n+1}:= Z\, C_n
- n\,\Theta \,(Y- \frac{n-1}{2}) C_{n-1} 
\end{equation}
and the series,
$$
\Phi(Z, \Theta)(s):= \sum \frac{s^n \, C_n}{n!}\, \Gamma(2Y+n)^{-1}
$$
\begin{lemma} $\Phi$ is the unique solution of the differential equation
$$
s (\frac{d}{ds})^2 \Phi - 2 (Y-1) \, \frac{d}{ds} \Phi 
+ Z \, \Phi - \frac{s}{2}\, \Theta \, \Phi=0
$$
which fulfills the further conditions,
$$
\Phi(0)=\,\Gamma(2Y)^{-1}
 \, , \quad  \frac{d}{ds} \Phi(0) = Z \,\Gamma(2Y+1)^{-1}
$$
\end{lemma}
\begin{proof} One has 
$$
 \frac{d}{ds} \Phi= \sum \frac{s^n \, C_{n+1}}{n!}\, \Gamma(2Y+n+1)^{-1}
$$
$$
 s (\frac{d}{ds})^2 \Phi= \sum \frac{s^n \, C_{n+1}}{n!}\, n \,\Gamma(2Y+n+1)^{-1}
$$
$$
2 (Y-1) \, \frac{d}{ds} \Phi=\sum \frac{s^n \, C_{n+1}}{n!}\,(2Y +2n)\, 
\Gamma(2Y+n+1)^{-1}
$$
$$
(n -(2Y +2n))\Gamma(2Y+n+1)^{-1}=-\Gamma(2Y+n)^{-1}
$$
so that 
$$
s (\frac{d}{ds})^2 \Phi - 2 (Y-1) \, \frac{d}{ds} \Phi 
=-\sum \frac{s^n \, C_{n+1}}{n!}\, \Gamma(2Y+n)^{-1}
$$
but by (\ref{indb1})
$$ 
C_{n+1}= Z C_n
- n\,\Theta  \,( Y- \frac{n-1}{2}) C_{n-1}
$$ 
the first term gives -$Z \, \Phi$ while the second gives,
$$
\Theta \,\sum \frac{s^n \, C_{n-1}}{n!}\,(n \, ( Y+ \frac{n-1}{2})) 
\Gamma(2Y+n)^{-1}
$$
which equals
$$
\frac{s \, \Theta}{2} \,\sum \frac{s^n \, C_n}{n!}\, ( 2Y+ n) 
\Gamma(2Y+n+1)^{-1}= \frac{s \, \Theta}{2} \,\Phi
$$

\end{proof}

\noindent  Let now $\mu$ be an operator in $\Ac$ such that 
\begin{equation} \label{X}
[ \Theta, \, \mu]= 0\, , \quad [Y, \, \mu]= \mu \, , \quad [[Z, \, \mu],\, \mu ]=0
\end{equation}
we let 
$$
\Psi(s) := e^{\frac{s\, \mu}{2}}\, \Phi(Z, \, \Theta)(s)
$$
\begin{lemma} $\Psi(s)$ satisfies the  following differential
equation
$$
s (\frac{d}{ds})^2 \Psi - 2 (Y-1) \, \frac{d}{ds} \Psi 
+ (Z + \mu \, Y)\, \Psi - \frac{s}{2}\, ( \Theta +[Z, \, \mu]
+\frac{\mu^2}{2}) \Psi=0
$$
\end{lemma}
\begin{proof} One has 
$$
 \frac{d}{ds} \Psi= \frac{\mu}{2}\Psi+ \,e^{\frac{s\, \mu}{2}}\,
\frac{d}{ds} \Phi
$$
$$
 s (\frac{d}{ds})^2 \Psi= s\, \frac{\mu^2}{4}\,\Psi
+ s \, \mu \, e^{\frac{s\, \mu}{2}}\,
\frac{d}{ds} \Phi+\, e^{\frac{s\, \mu}{2}}\,
 s (\frac{d}{ds})^2 \Phi
$$
$$
2 (Y-1) \, \frac{d}{ds} \Psi= \mu \, Y \, \Psi
+ s \, \mu \, e^{\frac{s\, \mu}{2}}\,
\frac{d}{ds} \Phi + \, e^{\frac{s\, \mu}{2}}\,
2 (Y-1) \, \frac{d}{ds} \Phi
$$
where for the last equality we used $[Y,\, \mu]=\, \mu$
to get
$$
[2Y,\, e^{\frac{s\, \mu}{2}}]=\,s\,\mu \, e^{\frac{s\, \mu}{2}}
$$
thus, 
$$
(s (\frac{d}{ds})^2-2 (Y-1) \, \frac{d}{ds})\Psi=
s\, \frac{\mu^2}{4}\,\Psi -\mu \, Y \, \Psi
+\, e^{\frac{s\, \mu}{2}}\,
(s (\frac{d}{ds})^2-2 (Y-1) \, \frac{d}{ds})\Phi
$$
but 
$$
(s (\frac{d}{ds})^2-2 (Y-1) \, \frac{d}{ds})\Phi
=-Z \, \Phi + \frac{s}{2}\,  \Theta \, \Phi
$$
and by (\ref{X}) one has,
$$
[Z, \, e^{\frac{s\, \mu}{2}}]=\frac{s}{2}\, [Z, \, \mu]
\, e^{\frac{s\, \mu}{2}}
$$
so that 
$$
 e^{\frac{s\, \mu}{2}}\,(-Z)= \,(-Z)\,e^{\frac{s\, \mu}{2}}
+\frac{s}{2}\, [Z, \, \mu]
\, e^{\frac{s\, \mu}{2}}
$$
and since $\mu$ commutes with $ \Theta$,
$$
(s (\frac{d}{ds})^2-2 (Y-1) \, \frac{d}{ds})\Psi=
s\, \frac{\mu^2}{4}\,\Psi -\mu \, Y \, \Psi
-Z \, \Psi +\frac{s}{2}\, [Z, \, \mu]\, \Psi
+ \frac{s}{2}\,  \Theta \, \Psi
$$
\end{proof}
Note also that
$$
\Psi(0)=\,\Gamma(2Y)^{-1}
 \, , \quad  \frac{d}{ds} \Psi(0) = (Z+ \, \mu \, Y) \,\Gamma(2Y+1)^{-1}
$$
It thus follows that the following holds.

\begin{proposition} \label{pert}
Let $\mu$ fulfill conditions  (\ref{X}), then 
$$ 
\Phi(Z+ \, \mu \, Y,\;  \Theta +\, [Z, \, \mu]\,+\, \frac{\mu^2}{2})(s)=\,  
e^{\frac{s\, \mu}{2}}\,\Phi(Z,  \Theta)(s)
 $$
\end{proposition}

\noindent One has by construction
$$
\Phi(X,  \Omega)(s)=\, \sum \frac{s^n \, B_n}{n!}\, \Gamma(2Y+n)^{-1}
$$
and similarly, 
$$
\Phi(S(X),  \Omega^{\,o})(s)=\, \sum \frac{s^n \, A_n}{n!}\, \Gamma(2Y+n)^{-1}
$$
where the $A_n$ are obtained by induction using,
\begin{equation} \label{inda1}
A_{n+1}:= S(X)\, A_n
- n\,\Omega^{\,o} \,(Y- \frac{n-1}{2}) A_{n-1} 
\end{equation}
while $A_{-1}:=0, \, A_{0}:=1$.

\noindent 
The first values of $A_n$ are
$$
A_1=  S(X)
$$
$$
A_2=  S(X)^2 - \,\Omega^{\,o} \, Y 
$$
$$
A_3=S(X)^3  - \,\Omega^{\,o} \,S(X) \,(3 Y+1) + \,X(\Omega)^{\,o} \, Y
$$
and 
\begin{eqnarray}
&& A_4=S(X)^4 
-\,\Omega^{\,o} \,S(X)^2 \, (6Y+4)+\,X(\Omega)^{\,o} \,S(X) \, (4Y+1)\nonumber \\
&&  \quad \quad \quad- \,X^2(\Omega)^{\,o} \, Y + 3\,(\Omega^{\,o})^2 \,Y \,
( Y+1)\nonumber
\end{eqnarray}

\noindent The general formula for $RC_n$ is
obtained as in (\ref{rc1}) by expanding the 
product
$$
\Phi(S(X), \Omega^{\,o})(s)(a)\, \Phi(X, \Omega)(s)(b)
$$
which gives,
\begin{equation}
RC_n(a,b):= \sum_{k=0}^n \, \frac{A_k }{k!} \,(2Y+k)_{n-k}(a)\,\,  
\frac{ B_{n-k}}{(n-k)!} \, (2Y+n-k)_{k}(b)   
\end{equation}

\begin{lemma} \label{mu1}
Let $\gamma \in Z^1(\Hc_1, \Ac)$ be 
a $1$-cocycle 
such that $\displaystyle
\gamma(X)=\gamma(Y)=0\,, \,\d_k(\gamma(h))=0\, \,\, \forall h\in \Hc_1\,,
\,\,k \in \mathbb{N}.
$
Then the brackets $RC_n$ are invariant under the inner
perturbation of the action of $\Hc_1$ associated to  $\gamma$.
\end{lemma}
\begin{proof}  Let $\displaystyle \mu :=\gamma(\d_1)$.
One has $\displaystyle \d_k(\mu)=0\, \,\, \forall k \in \mathbb{N}$.
The cocycle law
\begin{equation} \label{1coc}
\g (h \, h') = \sum
\g(h_{(1)}) \, h_{(2)}(\g(h')) \,,\quad \forall h \in \Hc \, .
\end{equation} 
shows that $\displaystyle \g(\d_2)= X(\mu) + \mu^2\,$
and $\,\displaystyle \g(\d'_2)=\, X(\mu) + \frac{1}{2}\mu^2\,$,
which using $\displaystyle [\d_1,\,\d'_2]=0\,$ gives
$\displaystyle [\mu,\,X(\mu)]=0$. We then get
\begin{equation} \label{cond}
[\Omega,\,\mu]=0\,, \quad Y(\mu) = \,\mu\,, \quad [\mu,\,X(\mu)]=0 
\end{equation}
 The effect of the perturbation on the generators is 
\begin{eqnarray}
&&  Y \rightarrow Y \,,\quad X \rightarrow X + L_\mu \, Y \,,\quad
\d_1 \rightarrow \d_1 + L_\mu - R_\mu\,,\nonumber \\
&&  \Omega  \rightarrow 
 \Omega +\,
X(\mu)\,+\, \frac{\mu^2}{2}
\end{eqnarray}
where $L_\mu$ is the 
operator  of left multiplication by $\mu$
and $R_\mu$ is right multiplication by $\mu$.

\noindent  We first apply Proposition \ref{pert} to the 
operator $L_\mu$ of left multiplication by $\mu$, and get,
using (\ref{cond}) to check (\ref{X}),
\begin{equation} \label{left}
\Phi(X + L_\mu \, Y,\; \Omega +\,
X(\mu)\,+\, \frac{\mu^2}{2}
)(s)= L_{e^{\frac{s\, \mu}{2}}}\, \Phi(X, \Omega)(s)
\end{equation}

\noindent The effect of the perturbation on $S(X)$ is 
\begin{equation}
S(X) \rightarrow \, S(X)-R_\mu \, Y,
\end{equation}
where $R_\mu$ is right multiplication by $\mu$.

\noindent One has
$$
\Phi(S(X), \Omega^{\, o})(s)= \sum \frac{s^n \, A_n}{n!}\, 
\Gamma(2Y+n)^{-1}
$$
\noindent  We now apply Proposition \ref{pert} to the 
operator $-R_\mu$ of right multiplication by $-\mu$,
using $\displaystyle \,S(X)(\mu)=\,-X (\mu)\,$ and (\ref{cond})
to check (\ref{X}) for the operators $\,\,\displaystyle S(X), \,Y, \, 
\Omega^{\,o}= R_\Omega,
\,-R_\mu$ and get,
\begin{equation} \label{right}
\Phi(S(X) -R_\mu \, Y,\; (\Omega +\,
X(\mu)\,+\, \frac{\mu^2}{2})^{\,o}
)(s)= R_{e^{-\frac{s\, \mu}{2}}}\, \Phi(S(X), \Omega)(s)
\end{equation}
combining (\ref{left}) and (\ref{right}) shows that the product
$$
\Phi(S(X), \Omega^{\,o})(s)(a)\, \Phi(X, \Omega)(s)(b)
$$
is unaltered by the perturbation and gives the required
invariance.
\end{proof}
\begin{lemma} \label{22} Let the Hopf algebra $\Hc_1$ act on an 
algebra $\Ac$ with $\d'_2$ inner as above
and $u \in \Ac$ be invertible and such that with $\mu =
u^{-1}\d_1(u)$,
$$
X(u)=0, \, Y(u)=0, \, \d_n(\mu)=0, \, \forall \, n \in \mathbb{N}
$$
\begin{itemize}
\item[$1^0$.] For all values of $n$,
$$
RC_n(x\,u,\,y)=\, RC_n(x,\,u\,y) \quad\forall x,\, y \in \Ac
$$

\item[$2^0$.] For all $x,\, y \in \Ac
$,
$$
RC_n(u\,x,\,y)=\,u\, RC_n(x,\,\,y)\, ,
\quad RC_n(x,\,y \,u)=\, RC_n(x,\,\,y)\, u  
$$

\item[$3^0$.] Let $\a$ be the inner automorphism implemented by $u$ one has,
$$
RC_n(\a(x),\,\a(y))=\, \a(RC_n(x,\,y)) \quad\forall x,\, y \in \Ac
$$
\end{itemize}
\end{lemma}
 \begin{proof} One has $Y(\mu)=\mu$, let us  show that
\begin{equation} \label{hyp}
[X(\mu),\,\mu] = 0
\end{equation}

 One has $\displaystyle \d_2(u)= [X,\,\d_1](u)=X(\d_1(u))
=X(u\, \mu)=u\,( X(\mu) + \mu^2)$, so that,
\begin{equation}  
\d'_2(u)=u\, \rho\,, \quad \rho := X(\mu) +\frac{1}{2} \mu^2
\end{equation}
Since $\displaystyle \d_1(X(\mu))=\, -[X,\,\d_1](\mu)=-\d_2(\mu)=0$, 
one has $\displaystyle \d_1(\rho)=0$.
The commutation  $\:\displaystyle [\d'_2 ,\, \d_1]\,(u)=0\,$ then entails
$$
u \, \rho\, \mu =\,  u \, \mu \, \rho \,,\quad [\mu, \, \rho]=0
$$
which implies (\ref{hyp}). Then as in (\ref{left1})
\begin{equation} \label{left10}
X\,L_u = L_u (\, X + L_{\mu}\, Y)
\end{equation}
where 
 $\mu :=u^{-1}\d_1(u)= -\d_1(u^{-1}) u$. 
Moreover,
\begin{equation} \label{left11}
\Omega\,L_u = L_u (\, \Omega + \, [X, L_{\mu}]\,+\frac{1}{2} L_{\mu}^2)
\end{equation}
which one gets using 
$$
\d'_2(u)= [\Omega, \, u]
$$
Since $\:\displaystyle [\Omega, \, \mu]= \, \d'_2(\mu) =0$
the hypothesis (\ref{X}) of Proposition \ref{pert} is fulfilled 
by $X,\,Y,\, \Omega,\,L_{\mu}\,$ and
we then get
$$ 
\Phi(X,\; \Omega)(s)L_u= \, L_u\,  L_{e^{\frac{s\, \mu}{2}}}\,
\Phi(X, \Omega)(s)
 $$
In a similar manner,
\begin{equation} \label{right10}
S(X)\,R_u = R_u (\, S(X) + R_{\mu'}\, Y)
\end{equation}
where 
 $\mu' :=\d_1(u)\,u^{-1}$. Moreover
\begin{equation} \label{right11}
\Omega^{\, o}\,R_u = R_u (\, \Omega^{\, o} + \, 
[S(X), R_{\mu'}]\,+\frac{1}{2}\,R_{\mu'}^2)
\end{equation}
where $\Omega^{\, o}$ really stands for $\,
R_{\Omega}$.
So  by Proposition \ref{pert}
we  get
$$ 
\Phi(S(X),\; \Omega^{\, o})(s))R_u= \, R_u\,  
R_{e^{\frac{s\, \mu'}{2}}}\,\Phi(S(X), \Omega)(s)\,.
 $$
Since $\;\displaystyle u\, \mu =\mu' \, u\;$
we conclude as in the proof of Lemma \ref{11}. 
\end{proof}

\noindent We then obtain the following naturality property of the 
construction of the higher brackets.

\bigskip
\begin{theorem} \label{cross}
When applied to any of the modular Hecke algebras 
$\Ac(\Gamma)$ the functor $RC_*$ yields the
reduced algebra of the crossed products 
of the algebra of modular forms endowed with the
Rankin-Cohen brackets by the action of the 
group $  \GL (2 , \mathbb{A}_{f})^0 $.
\end{theorem}
\begin{proof}
The algebra $\Ac(\Gamma)$ is the reduced algebra 
of $ \Mc \rtimes \GL (2 , \mathbb{A}_{f})^0 $
by the projection $e_\Gamma$ associated to $\Gamma$
(cf. \cite{CM10}).
With $u$ and $v$ as in Lemma \ref{22} and $a$ and $b \in \Ac$
one gets 
$$
RC_n(a\,u,\, b\,v)=\,RC_n(a,\,u\, b\,v)=\,RC_n(a,\,b^u) \,u\, \,v
$$
where $\displaystyle b^u:=u\,b\,u^{-1}$, which shows that
in the  crossed product algebra $\Ac= \Mc \ltimes 
\GL (2 , \mathbb{A}_{f})^0$ the $RC_n$ are entirely 
determined by their restriction to $\Mc$.
\end{proof}

\bigskip
\begin{corollary} \label{assoc}
When applied to any of the algebras 
$\Ac(\Gamma)$ the functor $RC_*$ yields associative 
deformations.
\end{corollary} 
\begin{proof}
The crossed product of an associative algebra by an 
automorphism group is associative, as well as its reduced algebras.
\end{proof}

\noindent Specifically, the product formula 
$$
a *_{t} b:= \sum {t}^n RC_n(a,b)
$$
gives an associative deformation.
More generally, following~\cite{CMZ}, 
for any $\kappa \in \mathbb{C}$ 
the product rule
\begin{equation} \label{fam}
a *_{t}^{\kappa} b \, := \,
\sum {t}^n \,  RC_n \left({\tb}_n^{\kappa}
(Y\ot 1, 1\ot Y)(a \ot b) \right) \, ,
\end{equation}
where  
\begin{equation} \label{fam1}
{\tb}_n^{\kappa}(x,y):= (-\frac{1}{4})^n\sum_j  \binom{n}{2j}
\frac{\binom{-\frac{1}{2}}{j}\binom{\kappa 
-\frac{3}{2}}{j}\binom{\frac{1}{2}
-\kappa}{j}}{\binom{-x-\frac{1}{2}}{j}\binom{-y-\frac{1}{2}}{j}
\binom{n+x +y -\frac{3}{2}}{j}}
\end{equation}
are the twisting coefficients defined in \cite{CMZ}, gives
an associative deformation.

\section{Rankin-Cohen Deformations}

We now return to the general case, we let the Hopf algebra 
$\Hc_1$ act on an algebra $\Ac$ and assume that
the derivation $\d'_2$ is inner implemented by 
an element $\Omega \in \Ac$,
\begin{equation}  
\d'_2(a)= \Omega \, a - a \, \Omega\,,\:\forall a \in \Ac
\end{equation}
with
\begin{equation}  
   \d_k(\Omega )=\,0\,,\: \forall k \in \mathbb{N}
\end{equation}
 Such an action of $\Hc_1$ on an algebra $\Ac$
will be said to define  a \textit{projective structure} on  $\Ac$,
and the element $\Om \in \Ac$ implementing the inner derivation
$\d'_2$ will be called its \textit{quadratic differential}.
\medskip

\noindent The main result of this section, extending Corollary \ref{assoc},
can be stated as follows.  
\medskip
 
\begin{theorem} \label{assoc1}
 The functor $RC_*$ applied to any algebra $\Ac$ endowed with
 a projective structure yields a  
 family of formal associative deformations of $\Ac$,  
 whose products are given by formula (\ref{fam}).
\end{theorem}
\bigskip

\noindent
In preparation for the proof, we shall extend the scalars
in the definition of $\Hc_1$.
Let $\Pc $ denote the free commutative algebra  
generated by the indeterminates 
$$ \{ Z_0, Z_1, Z_2, \ldots, Z_n,
\ldots \} .
$$
We define an action of  $ \Hc_1$ on $\Pc $
by setting on generators
\begin{equation} \label{Lie}
Y(Z_j) := (j+2)\,Z_j \, , \,  X(Z_j) := Z_{j+1} \, ,
  \quad \fl \, j \geq 0   \, ,
\end{equation}
and then extending $Y$, $X$ as derivations, while  
\begin{equation} 
\d_k (P) := 0 \, , \quad \fl \, P \in \Pc  \, .
\end{equation}
Equivalently, the Hopf action of $ \Hc_1$ is lifted from the
Lie algebra action defined by (\ref{Lie}).

\noindent We then form the double crossed product algebra
\begin{equation} 
\wt{\Hc}_1 = \Pc \rtimes \Hc_1 \ltimes \Pc \, ,
\end{equation}
whose underlying vector space is $\, \Pc \ot \Hc_1 \ot \Pc \,$,
with the product defined by the rule  
\begin{equation} \label{prod}
   P \rtimes h \ltimes Q \, \cdot \,  P' \rtimes h' \ltimes Q' 
   \, := \, \sum_{(h)} P \, h_{(1)} (P') \ltimes h_{(2)} h'
   \rtimes h_{(3)} (Q') \, Q \, ,
\end{equation}
where $\, P, Q, P', Q' \in \Pc \,$ and
$\, h, h' \in \Hc_{1} \,$ . 
\smallskip

\noindent
We next proceed to equip $\wt{\Hc}_1$ with the structure of an 
extended Hopf algebra over 
$\Pc $ (comp.~\cite{CM3}). First we turn $\wt{\Hc}_1$ into a (free)
$\Pc$-bimodule, by means of
the source and target homorphisms 
$\a \, , \, \b: \Pc \ra \wt{\Hc}_1$,  
\begin{equation} \label{st}
    \a (P) :=  P \rtimes 1  \ltimes 1 \, , \quad \text{resp.} \quad
    \b (Q) := 1 \rtimes 1 \ltimes Q \, , \quad \fl \, P, Q \in \Pc \, .
\end{equation}
 Note that
\begin{equation} \label{ab}
  P \rtimes h \ltimes Q \, = \, \a (P) \cdot \b (Q) \cdot h \, ,
  \qquad P, Q \in \Pc, \quad h \in \Hc_{1} \, .
\end{equation}
Note also that, while the commutation rules of $h \in \Hc_{1}$
with $\a (P)$ are given by the above action (\ref{Lie}) of 
$\Hc_{1}$ on $\Pc$, the commutation rules with $\b (Q)$
are more subtle ; for example (comp.~\cite[1.12]{CM3}),
\begin{equation} \label{Xb}
X \, \b (Q)\, - \, \b (Q)\, X \, = \, \b (X(Q))\,+ 
\, \b (Y(Q))\,\d_1 \, , \quad Q \in \Pc, \quad h \in \Hc_{1} \, .
\end{equation}
Let $\wt{\Hc}_1 \ot_{\Pc}\wt{\Hc}_1$ be the tensor square of
$\wt{\Hc}_1$ where we view $\wt{\Hc}_1$ as a bimodule over $\Pc$, 
using left multiplication by $ \b(\cdot)\,$
to define the right module structure and left multiplication by 
$ \a(\cdot)\,$ to define the left module structure.
\smallskip

\noindent The 
coproduct $\D : \wt{\Hc}_1 \ra \wt{\Hc}_1 \ot_{\Pc}\wt{\Hc}_1 $
is defined by
\begin{equation} \label{coprod}
    \D (P \rtimes h  \ltimes Q) := \sum_{(h)} P \rtimes h_{(1)} \ltimes
    1 \, \ot  \, 1 \rtimes h_{(2)} \ltimes Q 
\end{equation}
and satisfies the properties listed in~\cite[Prop. 6]{CM3}.
In particular, while the product of two elements in
$\wt{\Hc}_1 \ot_{\Pc} \wt{\Hc}_1 $ is not defined in general,
the fact that $\D$ is multiplicative, i.e. that
$$
\D (h_1 \cdot h_2 ) = \D (h_1) \cdot \D (h_2 ) \, , 
	\quad \quad \fl \, h_1 , h_2 \in \wt{\Hc}_1 \, ,
$$
makes perfect sense because of the property
$$ \D \, (h) 
 \cdot \left( \b (Q) \ot 1 \ - \ 1 \ot \a (Q) \right) = 0 \, ,
	\quad \fl \, Q \in \Pc \, , \,  h \in \wt{\Hc}_1 \, ,
$$
which uses only
the right action of $\wt{\Hc}_1 \ot \wt{\Hc}_1$ on
$\wt{\Hc}_1 \ot_{\Pc} \wt{\Hc}_1 $ by right multiplication.
In turn, since $\Hc_1$ is a Hopf algebra, it suffices to check
the latter on the algebra generators, i.e. for
$\, h = Y, X$ or $\d_1$. In that case
\begin{eqnarray*}
 \D (h) \cdot \left( \b (Q) \ot 1 \, - \, 1 \ot \a (Q) \right) \, = \nonumber \\
\sum_{(h)} \, \left(  h_{(1)} \cdot \b (Q) \ot_{\Pc} h_{(2)} \, - \,
h_{(1)} \cdot \ot_{\Pc}  h_{(2)} \a (Q) \right) \, = \nonumber \\
  \sum_{(h)} \, \left([ h_{(1)} , \b (Q) ] \ot_{\Pc} h_{(2)} \,
- \, h_{(1)} \ot_{\Pc}  [h_{(2)},\a (Q)] \right)  \, = \, 0 \, , \nonumber \\
\end{eqnarray*} 
where one needs (\ref{Xb}) to establish the vanishing.
\smallskip

\noindent 
The counit map $\ve: \wt{\Hc}_1 \ra \Pc$ is defined by
\begin{equation*}  
    \ve ( P \rtimes h \ltimes Q) := \, P \, \ve(h)\, Q \, , \quad
 P, Q \in \Pc, \quad h \in \Hc_{1} 
\end{equation*}
and fulfills the conditions listed in~\cite[Prop. 7]{CM3}.
\smallskip

\noindent Finally, the formula for the antipode is 
\begin{eqnarray*}  
S(P \rtimes h \ltimes Q) := S(h)_{(1)} (Q) \rtimes S(h)_{(2)} \ltimes
S(h)_{(3)} (P)  
 =   S(h)\cdot \a(Q) \cdot \b(P) .
\end{eqnarray*}
\medskip

\noindent In the same vein, an algebra $\Ac$
is a module-algebra over $\wt{\Hc}_1 \vert \Pc$ if,
first of all, $\Ac$ is gifted with an algebra homomorphism
$\rho: \Pc \ra \Ac$ (playing the role of the unit map over $\Pc$),
which turns $\Ac$ into a $\Pc$-bimodule via left and right
multiplication by the image of $\rho$, and secondly $\Ac$ is
endowed with an action  
$H \ot a \mapsto H(a) \, , \quad H \in \wt{\Hc}_1 \, , \,
a \in \Ac \,$ satisfying besides the usual action rules
\begin{eqnarray}
 (H\cdot H') (a) &=& H(H'(a)) \, , \qquad \qquad H, H' \in \wt{\Hc}_1
 \, , \  \\
   1 (a) &=& a \, , \qquad \qquad \qquad a \in \Ac \, , \nonumber 
\end{eqnarray}
also the compatibility rules
\begin{eqnarray} \label{comp}
 H(a_1 a_2) &=& \sum_{(H)} H_{(1)} (a_1) \, H_{(2)} (a_2) \, ,
 \qquad  a_1, a_2 \in \Ac \, , \\
 H (1) &=& \rho (\ve(H)) \, , \qquad \qquad H \in \wt{\Hc}_1 \, .
 \nonumber
\end{eqnarray} 
In particular for any $P \in \Pc$,
\begin{equation}
\a (P) (a) = \rho (P)\, a \, , \quad \text{resp.} \quad
\b (P) (a) = a \, \rho (P) \, , \qquad  \, a \in \Ac \, ,
\end{equation}
and therefore more generally for any monomial
$H = P \rtimes h \ltimes Q \in \wt{\Hc}_1$ one has
\begin{equation}
 P \rtimes h \ltimes Q \, (a) = \rho (P) \, h(a) \, \rho (Q) \, .
\end{equation}

\noindent We denote
\begin{equation} 
\wt{\d}'_2 \, := \, \d_2 - \frac{1}{2} \d_1^2 - \a(Z_0) + \b(Z_0) 
\end{equation}
and remark that it is a primitive element in $\wt{\Hc}_1$:
\begin{equation} \label{prim}
 \D (\wt{\d}'_2) \, = \, \wt{\d}'_2 \ot 1 \, + \, 1 \ot \wt{\d}'_2 \, .
\end{equation}

\noindent
We let $\wt{\Hc}_{\rm s}$ denote the quotient of $\wt{\Hc}_1$
by the ideal generated by $\wt{\d}'_2 $. In view of (\ref{prim}),
the latter is also a coideal, and therefore $\wt{\Hc}_{\rm s}$
inherits the structure of an extended Hopf algebra over $\Pc$.  
Clearly, the action of
$\wt{\Hc}_1 \vert \Pc$ on an algebra 
$\Ac$ endowed with a projective structure descends to
an action of $\wt{\Hc}_{\rm s}$ on $\Ac$.
\bigskip

\noindent As already mentioned in Section 3, the prototypical 
examples of projective structures are furnished
by the modular Hecke algebras of~\cite{CM10}.
For the purposes of this section it will suffice to consider
the `discrete' modular Hecke algebra, that is the crossed
product
$$  
\Ac_{G^+ (\mathbb{Q})} \, := \, \Mc \rtimes G^+(\mathbb{Q}) \, ,
\qquad G^+ (\mathbb{Q}) = \GL^+ (2, \mathbb{Q}) \, .
$$
where $\Mc$ is the algebra of modular forms of all levels.
We recall that
$\Ac_{G^+ (\mathbb{Q})} $ consists of finite sums of symbols of the form 
$$ \sum \, f \, U_{\g}^* \, , \qquad \text{with} \quad f \in \Mc \, ,
\quad \g \in  G^+ (\mathbb{Q})\, ,
$$
with the product given by the rule
\begin{eqnarray} 
f \, U_{\a}^* \,\cdot  g \, U_{\b}^*  =\,( f\,\cdot g| \a)
 \;  U_{\b \,  \a}^* \, ,
\end{eqnarray}
where the vertical bar denotes the `slash operation'. 
Under the customary identification
$$ f \in \Mc_{2k} \, \longmapsto \, \wt{f} := f \, dz^{k} 
$$
of modular forms with higher differentials, the latter is just
the pullback 
\begin{equation} \label{slash}
     f |\g \, \longmapsto  \, \g^{*} (\wt{f}) \, ,
    \quad \fl \, \g \in  G^+ (\mathbb{Q})\, .   
\end{equation}
We further recall that the
action of $\Hc_1$ on $\, \Ac_{G^+ (\mathbb{Q})} \,$ is determined by
\begin{equation} 
X (f \, U_{\g}^*) = X(f)\, U_{\g}^* \, , \,
Y (f \, U_{\g}^*) = Y(f)\, U_{\g}^* \, , \,
\d_1 (f \, U_{\g}^*) = \mu_{\g} \cdot f \, U_{\g}^* \, ,
\end{equation}
where
\begin{equation} \label{ram}
X(f) = \frac{1}{2  \pi i} \left(\frac{df}{dz} \, -
\, \frac{d}{dz} (\log \eta^4) \cdot Y(f) \right) \, ,
\end{equation}
$\, Y \,$ stands for the Euler operator, 
\begin{equation}
Y(f) = \frac{\ell}{2} \cdot f \, ,  
\end{equation}
for any $f \,$  of weight $\ell \,$. Lastly,
\begin{equation} \label{dmu}
\mu_{\g} \, (z) \, = \, 
\frac{1}{2  \pi i} \, \frac{d}{dz} \log \frac{\eta^4 | \g}{\eta^4} \, . 
\end{equation}
By \cite[Prop. 10]{CM10}, the quadratic differential is
implemented by the normalized 
Eisenstein modular form of level $1$ and weight $4$: 
\begin{equation} 
\d'_2 (a) = [\,\om, \, a]\, , \qquad  \om = \frac{E_4}{72}  \, .
\end{equation}
\medskip

\noindent
More generally, one can conjugate the above action
by a $G^+ (\mathbb{Q})$-invariant $1$-cocycle, 
as in~\cite[Prop. 11]{CM10}. In particular,
given any
$ \s \in \Mc_2 \,$,
there exists a unique $1$-cocycle
$\displaystyle u = u({\s}) \in 
\Lc( \Hc_1 ,\Ac_{G^+ (\mathbb{Q})})$ such that
\begin{equation} \label{cou}
u(X)= 0 , \;u(Y)= 0 , \; 
u(\delta_1)= \s  \, .
\end{equation}
The conjugate under $u({\s})$ of the
above action of $\Hc_1$ is given
on generators as follows:
\begin{eqnarray} \label{pa}
Y_{\s} = Y \, , \quad X_{\s} &=& X + \s \, Y  \, ,  \\
 (\delta_1)_{\s} (a) &=& \, \delta_1 (a)
+ [{\s}, \, a] \, , \qquad a \in \Ac_{G^+ (\mathbb{Q})}\,. 
\nonumber 
\end{eqnarray}
The conjugate under $u({\s})$  of $\delta'_2$ is given
by the operator
\begin{equation} \label{us}
(\delta'_2)_{\s} (a) = [\,\om_{\s}, \, a] \quad \text{with} \quad
\om_{\s} = \om + X(\s) + \frac{{\s}^2}{2}  \, ,
\quad a \in \Ac_{G^+ (\mathbb{Q})}.
\end{equation}
\medskip

\noindent Thus, for any $1$-cocycle $u = u({\s})$ as above,
we get an action of $\wt{\Hc}_{\rm s} \vert \Pc$ on 
$\Ac_{G^+ (\mathbb{Q})}$, determined by (\ref{pa})
and by the homomorphism  $\rho_{\s} : \Pc \ra \Mc$,
\begin{equation} \label{ro}
 \rho_{\s} (Z_k) \, = \, X_{\s}^k (\om_{\s} ) \, , \qquad
 k=0, 1, 2, \ldots \, .
\end{equation}
\medskip

\noindent For each $n$ we let $\rho_{\s}^{\ot n}
: \Pc^{\ot n}  \ra \Mc^{\ot n}$ 
be the nth tensor power of $\rho_{\s}$.
We shall first show that the  family $\rho_{\s}$
is sufficiently large to separate the elements of $\Pc^{\ot n}$.

\begin{lemma} \label{tech0}
For each $n \in \mathbb{N}$, one has
$\displaystyle \bigcap_{\s \in \Mc_2}  \Ker \rho_{\s}^{\ot n} = 0 \, $. 
 \end{lemma}
\medskip

\begin{proof} Let us treat the case $n=1$ first. 
Let $g_2^*$ be the ``quasimodular'' (\cite{Zagier}) 
   solution of the equation
\begin{equation} \label{zz}
X(m) +\frac{ m^2}{2}+\, \om \,=0 .
\end{equation}
Note that not only $g_2^*$ fulfills
(\ref{zz}) but also one has
$$
X(f) =\frac{1}{2  \pi i} \frac{d}{dz}\, f -g_2^* \, Y\, f \, ,
$$
for any modular form $f$. Also
$$
\frac{1}{2  \pi i} \frac{d}{dz}\,g_2^*- \frac{1}{2} (g_2^*)^2 =- \, \om
$$
 Thus, with $\a := \s -g_2^*$ we get 
\begin{eqnarray} \label{pa1}
  X_{\s} &=& \frac{1}{2  \pi i} \frac{d}{dz}\, + \a \, Y  \, ,  \\
\om_{\s} &=&  \frac{1}{2  \pi i} \frac{d }{dz}\, \a + \frac{{\a}^2}{2} 
\nonumber 
\end{eqnarray}
which allows to rewrite (\ref{ro}) as
\begin{equation} \label{ro1}
 \rho_{\s} (Z_k) \, = \, (\frac{1}{2  \pi i} 
 \frac{d}{dz} + \a \, Y)^k (\frac{1}{2  \pi i} \frac{d }{dz}\, 
 \a + \frac{{\a}^2}{2} 
 ) \, , \qquad
 k=0, 1, 2, \ldots \, .
\end{equation}
Given $0 \neq P \in \Pc$ the set of $\a$ for which $\rho_{\s}(P)=0$
is seen using (\ref{ro1}) to be contained in the space of holomorphic
solutions of an (autonomous) ODE \footnote{Any given quasi-modular form $\a$ is
the solution of a non-trivial autonomous ODE (cf. \cite{Zagier}),
but the latter of course depends on $\a$.} and these only depend
on finitely many parameters. Thus, given $0 \neq P \in \Pc$,
the set of $\s$ for which $\rho_{\s}(P)=0$ 
is finite dimensional.
\smallskip
 
\noindent Let us now prove by induction that
the same result holds for any $n$.
Given $P \in \Pc^{\ot n}$ we write 
$$
P = \sum P_j \ot m_j
$$
where $P_j 
\in 
\Pc^{\ot n-1}$ and the $
m_j$ belong to the canonical basis 
of monomials in $\Pc$.
If $P \neq 0$ then $P_j
\neq 0$ for some $j$ and by the induction hypothesis 
the set $E_j$
of $\s$ for which $\rho_{\s}(P_j)=0$
is finite dimensional. On the complement of $E_j$
any $\a := \s -g_2^*$ such that $\rho_{\s}(P)=0$ 
fulfills a non-trivial (autonomous) ODE
of the form 
\begin{equation} \label{ode}
\sum \lambda_j
 \rho_{\s} (m_j) \, =0
\end{equation}
where the coefficient $\lambda_j \neq 0$.
Since the space of parameters for equations of the form (\ref{ode})
is finite dimensional we get the required finite dimensionality. 
Since the space of modular forms of weight $2$ and arbitrary level
is infinite dimensional we conclude the proof.        
\end{proof}

\noindent We next define a map of bimodules
$$ \chi_{\s}^{(n)}: \, \underbrace{\wt{\Hc}_{\rm s} \ot_{\Pc} \ldots 
\ot_{\Pc}  \wt{\Hc}_{\rm s}}_{n-\text{times}} \quad 
\longrightarrow  \quad 
\Lc ( \underbrace{\Ac_{G^+ (\mathbb{Q})} \ot \ldots 
\ot \Ac_{G^+ (\mathbb{Q})}}_{n-\text{times}}  \, ,
\Ac_{G^+ (\mathbb{Q})} )
$$
by means of the assignment
\begin{equation} \label{chin}
\chi_{\s}^{(n)} (h^{1} \ot_{\Pc} \ldots \ot_{\Pc} h^{n}) \, 
(a_{1}, \ldots , a_{n}) 
\, = \,
 \  h_{\s}^{1}(a_{1}) \cdots h_{\s}^{n}(a_{n}) \, ,
\end{equation}
where $\, a_{1}, \ldots , a_{n} \in \Ac_{G^+ (\mathbb{Q})} \,$
and with $\, h^{1}, \ldots , h^{n} \in \wt{\Hc}_{\rm s}\, $ 
acting via $\rho_{\s}$.
\bigskip

\noindent The following result, which represents 
a `modular' analogue of
\cite[Prop.4]{CM3}, allows to
establish the associativity at the Hopf algebraic level.
\medskip

\begin{proposition} \label{tech}
For each $n \in \mathbb{N}$, one has
$$\bigcap_{\s}\,  \Ker \chi_{\s}^{(n)} = 0 \, . 
$$
\end{proposition}
\medskip

\begin{proof}  
For the sake of clarity, we shall first treat the
case $n=1$. 
An arbitrary element of $\wt{\Hc}_{\rm s}$ can be represented
uniquely as a finite sum of the form 
$$ H \, = \, \sum_{j, k, m, s \geq 0} \a(P_{jkms}) \, 
\b(Q_{jkms}) \, \d_1^j \, X^k \, Y^m \, ,
$$
with $P, Q \in \Pc$. 
Let us assume $\chi_{\s}^{(1)} (H) = 0$, for any $\s$.
\smallskip

\noindent
Evaluating $H$ on a generic monomial in $\Ac_{G^+ (\mathbb{Q})}$ 
one obtains, for any $f \in \Mc$ and any $\g \in G^+ (\mathbb{Q})$,
\begin{equation} \label{nul1}
\sum_{j, k, m, s \geq 0} P_{jkms} \, 
Q_{jkms}\vert \g \cdot \mu_{\g}(\s)^j \cdot
X^k (Y^m (f)) \, = \, 0 \, ;
\end{equation} 
where, 
\begin{equation}
 \mu_{\g}(\s) \, (z) = 
   \a ( z) - \a | \g \, ( z)  
- \, \frac{1}{\pi i}\,\,\frac{ c}{cz+d} 
\end{equation} 
with $\a := \s -g_2^*$ as above.
By continuity, the above  
holds in fact for any $\g \in G^+ (\mathbb{R})$.

\noindent
For each fixed $l$, the differential equation
\begin{equation*} 
\sum_{j, k,m, s} P_{jkms} 
Q_{jkms}\vert \g \, \mu_{\g}(\s)^j \, l^m \, X^k (f) 
\, = \, 0 
\end{equation*}
is satisfied by all modular forms $f \in \Mc_{2l}$. 
In turn, this implies that
all its coefficients vanish. Using the freedom in $l$, 
it then follows that  
\begin{equation} \label{nul2}
\sum_{j, s} P_{jkms} 
Q_{jkms}\vert \g \, \mu_{\g}(\s)^j
\, = \, 0 \, ,
\end{equation}
for each $k$ and $m$. 
\smallskip
 
\noindent 
Given $z \in H$, the following
three functions on $\SL(2,\mathbb{C})$ are defined and independent
in a neighborhood of $\SL(2,\mathbb{R})$:
\begin{eqnarray} 
g_1(a,b,c,d):&=& \frac{az+b}{cz+d} \, , \quad
g_2(a,b,c,d):= \frac{1}{cz+d} \, ,  \nonumber\cr \cr
g_3(a,b,c,d):&=& \frac{c}{cz+d} \,  . \nonumber
\end{eqnarray}
Indeed, we recover $c$ from $g_2$ and $g_3$
and then $d$ from $g_3$, then $az+b$ from
$g_1$ and finally $a$ and $b$ from $ad-bc=1$.
Now the formula for $\mu_{\g}$ is of the
form
\begin{equation}
\mu_{\gamma}(\s)= \a(g_1) (g_2)^2 - \a(z) -\frac{1}{i \,   \pi} g_3
\end{equation} 
and thus involves $g_3$ nontrivially while
the other terms in the formulas do not involve $g_3$.
For fixed $z$ the formula (\ref{nul2}) remains valid
in a neighborhood of $\SL(2,\mathbb{R})$ in $\SL(2,\mathbb{C})$.
Fixing $z \, , g_1 \, , g_2$ and varying $g_3$ independently
this is enough to show that the coefficient of each power
$\mu_{\g}(\s)^j$ vanishes identically. Thus, 
\begin{equation} \label{Null}
\sum_{s} P_{jkms} \cdot
Q_{jkms} \vert \g  \, = \, 0 \, , \qquad \fl \, 
\g \in \SL(2,\mathbb{R}) \, .
\end{equation}
Using the independence of the functions $\, g_1 \, , g_2 \,$
and the fact that the identity (\ref{Null}) holds true for every
homogeneous component of $Q_{jkms}$, in view of the freedom
to choose $\s$, it follows from Lemma \ref{tech0} that $H=0$.
\medskip

\noindent The proof of the case $n > 1$ is obtained by combining
the above arguments with the proof of the general case in
Lemma \ref{tech0}.
An arbitrary element of $\wt{\Hc}_{\rm s}^{\ot n}$ can be represented
uniquely as a finite sum of the form 
\begin{eqnarray} 
 H \, = \, \sum_{j, k, m, s } \a(P_{1,j_1k_1m_1s}) \, 
\b(Q_{1,j_1k_1m_1s}) \, \d_1^{j_1} \, X^{k_1} \, Y^{m_1} \,\ot \cdots  
\nonumber\cr \cr
\cdots \ot \, \a(P_{a,j_ak_am_as}) \, 
\b(Q_{a,j_ak_am_as}) \, \d_1^{j_a} \, X^{k_a} \, Y^{m_a}
\,\ot \cdots  \nonumber\cr \cr
\cdots
 \ot \, \a(P_{n,j_nk_nm_ns}) \, 
\b(Q_{n,j_nk_nm_ns}) \, \d_1^{j_n} \, X^{k_n} \, Y^{m_n}
\end{eqnarray}
\noindent
Evaluating $H$ on a generic monomial in $\Ac_{G^+ (\mathbb{Q})}^{\ot n}$ 
one obtains, for any $f_j \in \Mc$ and any $\g_j \in G^+ (\mathbb{Q})$,
\begin{equation}  
\sum_{j, k, m, s } \prod_a \,(P_{a,j_ak_am_as} \,\cdot
 Q_{a,j_ak_am_as}\vert \g_a \cdot
 \mu_{\g_a}(\s)^{j_a} \cdot
X_{\s}^{k_a} (Y^{m_a} (f_a)))\vert\g_{a-1} \cdots \g_1 \, = \, 0 \, 
\end{equation} 
As in the case $\, n=1$, the freedom in the choice of the $f_j \in \Mc$
and $\g_j \in G^+ (\mathbb{Q})$ shows that for every multiindex
$(j_a,k_a,m_a)_{a \in \{1,..,n\}}$, one has 
\begin{equation}  
\sum_{ s } \prod_a \,(P_{a,j_ak_am_as} \,\cdot
 Q_{a,j_ak_am_as}\vert \g_a )\vert\g_{a-1} \cdots \g_1 \, = \, 0 \, .
\end{equation} 
Using once more the independence of the 
functions $g_1$ and $g_2$ together with the reduction to homogeneous 
components for $ Q_{a,j_ak_am_as}$, and then applying
again Lemma \ref{tech0},
one arrives at the conclusion that $H=0$.
\end{proof}
\bigskip

\noindent As a direct consequence, one obtains a
family of `universal deformation formulas' (in the sense of
\cite{GZ}) based on $\wt{\Hc}_{\rm s}$, as follows.
From Theorem \ref{cross} and Proposition \ref{tech}, each
bilinear operator $RC_n$ uniquely determines
an element
\begin{equation}
 RC^{\Hc_{\rm s}}_n  \in 
\wt{\Hc}_{\rm s}   \ot_{\Pc } \wt{\Hc}_{\rm s}   \, .
\end{equation}
We then assemble all of them into
the element
\begin{equation}
\Fc \, := \,
\sum_{n \geq 0} \, {t}^n  RC^{\Hc_{\rm s}}_n  \, \in 
\wt{\Hc}_{\rm s} [[t]] \ot_{\Pc [[t]]} \wt{\Hc}_{\rm s} [[t]]
\end{equation}
\begin{corollary}\label{check} The element
$\Fc \in \wt{\Hc}_{\rm s} [[t]] \ot_{\Pc [[t]]} \wt{\Hc}_{\rm s} [[t]] $ 
defines a
universal deformation formula, i.e.  
satisfies the identities 
\begin{eqnarray} \label{twist}
  (\D \ot \Id) (\Fc) \cdot (\Fc \ot 1) \, &=& \, 
      (\Id \ot \D) (\Fc) \cdot 1 \ot \Fc  \, ; \\ 
  (\ve \ot \Id) (\Fc) \, &=& \,  1 \ot 1 \, = \, (\Id \ot \ve) (\Fc) \, .    
\end{eqnarray}
\end{corollary}
\medskip

\begin{proof} Indeed, this follows by applying
Proposition \ref{tech}, for $n=3$, to the  
associative deformation 
of the algebra $\Ac_{G^+ (\mathbb{Q})}$
provided by Corollary \ref{assoc}. 
\end{proof}
\medskip

\noindent
This result, in conjunction with  
the twisting operators ${\tb}_n^{\kappa}$ of (\ref{fam1}), 
immediately imply Theorem 10.
\medskip

\noindent \textbf{Remark \, 14.} We finally note that by evaluating $\Fc$
at $Z_k = 0 \, , k \geq 0 \,$, 
one obtains a completely explicit twisting element
$F \in \Hc_{\rm s} [[t]] \ot_{\mathbb{C} [[t]]} \Hc_{\rm s} [[t]] $,
\begin{equation}
 F \, = \, \sum_{n\geq 0} \, t^n \, \sum_{k=0}^n \, \frac{S(X)^k}{k!} \,
(2Y+k)_{n-k} \ot \frac{X^{n-k}}{(n-k)!} \,(2Y+n-k)_{k} \, ,  
\end{equation}
which defines a deformation  
of the Hopf algebra $\Hc_{\rm s}$
in the direction of the Hochschild $2$-cocycle
$-2T = -X \ot 2Y + 2Y \ot X + \d_1 \cdot Y \ot 2Y \,$.
Although the expressions of its components were found early on,
cf. (\ref{RC0}), the proof of their universal associativity 
property passes through
the treatment of the more general case
of an arbitrary quadratic differential.

\section{Appendix: Explicit Formulas}

\noindent We display below the formulas of $RC_n$ for the
first three values of $\, n$, 
illustrating the rapid increase of their complexity.
They are in reduced form, with the terms in $\a[\cdot]$ followed
by those in $\b[\cdot]$ appearing in the role of coefficients,
and then followed by
terms in $\d_1$, in $X$ and finally in $Y$; the last three types 
form the analogue of a Poincar\'e-Birkhoff-Witt basis of 
$\wt{\Hc}_{\rm s}$,
viewed as a module over $\Pc \ot \Pc$.
Moreover,
taking advantage of the tensoring over $\Pc$, no term in $\b[\cdot]$
appears in the first argument of the tensor square. 
\bigskip

\noindent  In the formulas that follow we shall lighten
the notation, using the symbol $RC_*$ instead of 
$RC_*^{\Hc_{\rm s}}$, and $\ot$ instead of $\ot_{\Pc}$; we
shall also write $\Om$ instead of $Z_0$.
 \bigskip
 
\begin{eqnarray} 
\frac{1}{2} RC_1 \, = \, -X \ot Y + Y \ot X + \d_1 \cdot Y \ot Y \, \, . 
\end{eqnarray}

\noindent Quite remarkably, the formula for $A_n$
appears to contain no term in $\b[.]$, unlike the formula
for $S(X)^n$. The latter, once put in reduced form, is much
more complicated than the expression of $A_n$.
\medskip

\noindent The associativity for $RC_1$ follows 
directly from its Hochschild property. It is already harder to 
check it directly for $RC_2$
which is given by the following expression: 
\medskip

\begin{eqnarray*} 
& & RC_2 \,= \, -X\otimes X-2\, X\otimes X\cdot Y+{X^2}\otimes Y+2\, 
{X^2}\otimes {Y^2}+  Y\otimes {X^2} 
\cr \cr & &
-\,Y\otimes \alpha [\Omega ]\cdot Y+2\, {Y^2}\otimes {X^2}-2\, 
 {Y^2}\otimes \alpha [\Omega ]\cdot Y- 2\, X\cdot Y\otimes X 
\cr \cr & &
 -4\, X\cdot Y\otimes
 X\cdot Y-{{\delta }_1}\cdot X\otimes Y-2\, {{\delta }_1}\cdot X\otimes {Y^2}+
  {{\delta }_1}\cdot Y\otimes X \cr \cr & &
+\, 2\, {{\delta }_1}\cdot Y\otimes X\cdot Y+2\, 
 {{\delta }_1}\cdot {Y^2}\otimes X+  4\, {{\delta }_1}\cdot {Y^2}\otimes
 X\cdot Y + \frac{1}{2}\, \delta _{1}^{2}\cdot Y\otimes Y \cr \cr & &
+\,\,\delta _{1}^{2}\cdot
Y\otimes {Y^2}+  
 \delta _{1}^{2}\cdot {Y^2}\otimes Y+2\, \delta _{1}^{2}\cdot {Y^2}\otimes 
 {Y^2}-\alpha [\Omega ]\cdot Y\otimes Y \cr \cr & &
- \,2\, \alpha [\Omega ]\cdot Y\otimes {Y^2}-2\, 
 {{\delta }_1}\cdot X\cdot Y\otimes Y-4\,
 {{\delta }_1}\cdot X\cdot Y\otimes {Y^2} \, .
\end{eqnarray*} 
\medskip

\normalsize
\noindent We did check directly Corollary \ref{check} up to order $4$
included, (i.e. the associativity for $RC_3$ and $RC_4$) 
with the help of a computer. This is 
beyond the reach of any `bare hands' computation, 
as witnessed by the complexity
of the following
formula for $RC_3$. (The expression of  $RC_4$ is much longer,
it would occupy several pages.)
 
\scriptsize
 \begin{eqnarray*} 
& & RC_3 \,
=-2\, X\otimes {X^2}-2\, X\otimes {X^2}.Y+2\, X\otimes \alpha [\Omega ].Y+2\, X\otimes \alpha [\Omega ].{Y^2}+  
 2\, {X^2}\otimes X+6\, {X^2}\otimes X.Y  +  \cr \cr & &
 + \,4\, {X^2}\otimes X.{Y^2}-\frac{2\, {X^3}\otimes Y}{3} -   
 2\, {X^3}\otimes {Y^2}-\frac{4\, {X^3}\otimes {Y^3}}{3}+\frac{2\, Y\otimes {X^3}}{3}-\frac{2}{3}\, Y\otimes \alpha
[\Omega ].X  \cr \cr & & -  
 \frac{2}{3}\, Y\otimes \alpha [X[\Omega ]].Y-2\, Y\otimes \alpha [\Omega ].X.Y
+  
 2\, {Y^2}\otimes {X^3}-2\, {Y^2}\otimes \alpha [\Omega ].X-2\, 
 {Y^2}\otimes \alpha [X[\Omega ]].Y  \cr \cr & &-  
 6\, {Y^2}\otimes \alpha [\Omega ].X.Y
+\frac{4\, {Y^3}\otimes {X^3}}{3}-\frac{4}{3}\, {Y^3}\otimes \alpha [\Omega ].X- 
 \frac{4}{3}\, {Y^3}\otimes \alpha [X[\Omega ]].Y-4\, {Y^3}\otimes \alpha [\Omega ].X.Y   \cr \cr & &-6\, X.Y\otimes {X^2}-  
 6\, X.Y\otimes {X^2}.Y+6\, X.Y\otimes \alpha [\Omega ].Y+6\, X.Y\otimes \alpha [\Omega ].{Y^2}-  
 4\, X.{Y^2}\otimes {X^2} \cr \cr & &-4\, X.{Y^2}\otimes {X^2}.Y  + \,4\, X.{Y^2}\otimes \alpha [\Omega ].Y+  
 4\, X.{Y^2}\otimes \alpha [\Omega ].{Y^2}+2\, {X^2}.Y\otimes X+6\, {X^2}.Y\otimes X.Y  \cr \cr & &+  
 4\, {X^2}.Y\otimes X.{Y^2}-2\, {{\delta }_1}.X\otimes X -6\, {{\delta }_1}.X\otimes X.Y-  
 4\, {{\delta }_1}.X\otimes X.{Y^2}+2\, {{\delta }_1}.{X^2}\otimes Y+6\, {{\delta }_1}.{X^2}\otimes {Y^2}\cr \cr & &+  
 4\, {{\delta }_1}.{X^2}\otimes {Y^3}+2\, {{\delta }_1}.Y\otimes {X^2}   + \,2\, {{\delta }_1}.Y\otimes {X^2}.Y-  
 2\, {{\delta }_1}.Y\otimes \alpha [\Omega ].Y-2\, {{\delta }_1}.Y\otimes \alpha [\Omega ].{Y^2}\cr \cr & &+  
 6\, {{\delta }_1}.{Y^2}\otimes {X^2}+6\, {{\delta }_1}.{Y^2}\otimes {X^2}.Y   -6\, {{\delta }_1}.{Y^2}\otimes \alpha [\Omega
].Y-  
 6\, {{\delta }_1}.{Y^2}\otimes \alpha [\Omega ].{Y^2}+4\, {{\delta }_1}.{Y^3}\otimes {X^2}\cr \cr & &+4\, {{\delta }_1}.{Y^3}\otimes
{X^2}.Y-  
 4\, {{\delta }_1}.{Y^3}\otimes \alpha [\Omega ].Y   -4\, {{\delta }_1}.{Y^3}\otimes \alpha [\Omega ].{Y^2}-  
 \delta _{1}^{2}.X\otimes Y-3\, \delta _{1}^{2}.X\otimes {Y^2} \cr \cr & &-2\, \delta _{1}^{2}.X\otimes {Y^3}+\delta _{1}^{2}.Y\otimes
X+  
 3\, \delta _{1}^{2}.Y\otimes X.Y  + \,2\, \delta _{1}^{2}.Y\otimes X.{Y^2}+3\, \delta _{1}^{2}.{Y^2}\otimes X+  
 9\, \delta _{1}^{2}.{Y^2}\otimes X.Y\cr \cr & &+6\, \delta _{1}^{2}.{Y^2}\otimes X.{Y^2}+2\, \delta _{1}^{2}.{Y^3}\otimes X+  
 6\, \delta _{1}^{2}.{Y^3}\otimes X.Y   + \,4\, \delta _{1}^{2}.{Y^3}\otimes X.{Y^2}+\frac{1}{3}\, \delta _{1}^{3}.Y\otimes
Y+  
 \delta _{1}^{3}.Y\otimes {Y^2} \cr \cr & &+\frac{2}{3}\, \delta _{1}^{3}.Y\otimes {Y^3}+\delta _{1}^{3}.{Y^2}\otimes Y+3\, \delta _{1}^{3}.{Y^2}\otimes
{Y^2}+  
 2\, \delta _{1}^{3}.{Y^2}\otimes {Y^3}  +\frac{2}{3}\, \delta _{1}^{3}.{Y^3}\otimes Y+2\, \delta _{1}^{3}.{Y^3}\otimes
{Y^2} \cr \cr & &+  
 \frac{4}{3}\, \delta _{1}^{3}.{Y^3}\otimes {Y^3}+ \,\frac{2}{3}\, \alpha [\Omega ].X\otimes Y+2\, \alpha [\Omega ].X\otimes
{Y^2}+  
 \frac{4}{3}\, \alpha [\Omega ].X\otimes {Y^3}   -2\, \alpha [\Omega ].Y\otimes X  \cr \cr & &-6\, \alpha [\Omega ].Y\otimes X.Y-  
 4\, \alpha [\Omega ].Y\otimes X.{Y^2}-2\, \alpha [\Omega ].{Y^2}\otimes X-6\, \alpha [\Omega ].{Y^2}\otimes X.Y-  
 4\, \alpha [\Omega ].{Y^2}\otimes X.{Y^2}   \cr \cr & &+\frac{2}{3}\, \alpha [X[\Omega ]].Y\otimes Y+  
 2\, \alpha [X[\Omega ]].Y\otimes {Y^2}+\frac{4}{3}\, \alpha [X[\Omega ]].Y\otimes {Y^3}-  
 6\, {{\delta }_1}.X.Y\otimes X-18\, {{\delta }_1}.X.Y\otimes X.Y   \cr \cr & &-12\, {{\delta }_1}.X.Y\otimes X.{Y^2}-  
 4\, {{\delta }_1}.X.{Y^2}\otimes X-12\, {{\delta }_1}.X.{Y^2}\otimes X.Y-  
 8\, {{\delta }_1}.X.{Y^2}\otimes X.{Y^2}+2\, {{\delta }_1}.{X^2}.Y\otimes Y   \cr \cr & &+   \,
 6\, {{\delta }_1}.{X^2}.Y\otimes {Y^2}+4\, {{\delta }_1}.{X^2}.Y\otimes {Y^3}-3\, \delta _{1}^{2}.X.Y\otimes Y-  
 9\, \delta _{1}^{2}.X.Y\otimes {Y^2}-6\, \delta _{1}^{2}.X.Y\otimes {Y^3}\cr \cr & &-2\, \delta _{1}^{2}.X.{Y^2}\otimes Y   -  
 6\, \delta _{1}^{2}.X.{Y^2}\otimes {Y^2}-4\, \delta _{1}^{2}.X.{Y^2}\otimes {Y^3}+2\, \alpha [\Omega ].X.Y\otimes Y+ 
 6\, \alpha [\Omega ].X.Y\otimes {Y^2}\cr \cr & &+4\, \alpha [\Omega ].X.Y\otimes {Y^3}   -  
 2\, \alpha [\Omega ].{{\delta }_1}.Y\otimes Y-6\, \alpha [\Omega ].{{\delta }_1}.Y\otimes {Y^2}-  
 4\, \alpha [\Omega ].{{\delta }_1}.Y\otimes {Y^3}-2\, \alpha [\Omega ].{{\delta }_1}.{Y^2}\otimes Y \cr \cr & &-  
 6\, \alpha [\Omega ].{{\delta }_1}.{Y^2}\otimes {Y^2}  -4\, 
 \alpha [\Omega ].{{\delta }_1}.{Y^2}\otimes {Y^3} \, \, .
\end{eqnarray*} 

\end{document}